\documentclass[12pt]{article}
\usepackage{fullpage,graphicx,psfrag,amsmath,amsfonts,verbatim}
\usepackage[small,bf]{caption}
\usepackage{amssymb}
\graphicspath{{figures/}}

\usepackage{listings}
\usepackage{longtable}
\usepackage{subcaption}

\usepackage{color}

\definecolor{codegreen}{rgb}{0,0.6,0}
\definecolor{codegray}{rgb}{0.5,0.5,0.5}
\definecolor{codepurple}{rgb}{0.58,0,0.82}
\definecolor{backcolour}{rgb}{0.95,0.95,0.98}

\lstdefinestyle{mystyle}{
backgroundcolor=\color{backcolour},
keywordstyle=\color{magenta},
numberstyle=\tiny\color{codegray},
stringstyle=\color{codepurple},
basicstyle=\ttfamily\small,
breakatwhitespace=false,
breaklines=true,
captionpos=t,
keepspaces=true,
numbers=left,
numbersep=5pt,
showspaces=false,
showstringspaces=false,
showtabs=false,
tabsize=2,
}

\definecolor{seagreen}{rgb}{0.18, 0.55, 0.34}
\definecolor{mediumviolet-red}{rgb}{0.78, 0.08, 0.52}
\definecolor{khaki}{rgb}{0.94, 0.9, 0.55}

\lstdefinelanguage{mypython}
{
keywords=[1]{from, import, as, assert, not, print, nonneg, PSD, axis},
keywordstyle=[1]{\color{mediumviolet-red}},
keywords=[2]{cp, lo, pl, cvxpy, Variable, Parameter,
sqrt, exp, numpy, np, Problem, Minimize, Maximize, value, solve, inner,
sum, multiply, arange, range, norm1, norm2, norm_inf, abs, square,
diagonal, outer, pos, hstack, power},
keywordstyle=[2]{\color{seagreen}},
upquote=true,
showstringspaces=false,
basicstyle=\ttfamily,
columns=fullflexible,
keepspaces=true,
emph={True,False,def,return,float,class,match,switch,len},
emphstyle={\color{seagreen}},
belowskip=1em,
aboveskip=1em,
morecomment=[l]{\#}
}

\newcommand{\BEAS}{\begin{eqnarray*}}
\newcommand{\EEAS}{\end{eqnarray*}}
\newcommand{\BEA}{\begin{eqnarray}}
\newcommand{\EEA}{\end{eqnarray}}
\newcommand{\BEQ}{\begin{equation}}
\newcommand{\EEQ}{\end{equation}}
\newcommand{\BIT}{\begin{itemize}}
\newcommand{\EIT}{\end{itemize}}
\newcommand{\BNUM}{\begin{enumerate}}
\newcommand{\ENUM}{\end{enumerate}}

\newcommand{\BA}{\begin{array}}
\newcommand{\EA}{\end{array}}

\newcommand{\ie}{{\it i.e.}}

\newcommand{\ones}{\mathbf 1}

\newcommand{\reals}{{\mbox{\bf R}}}

\newcommand{\symm}{{\mbox{\bf S}}}  % symmetric matrices

{\begin{quote}}{\end{quote}}

\makeatletter
\long\def\@makecaption#1#2{
\vskip 9pt
\begin{small}
\setbox\@tempboxa\hbox{{\bf #1:} #2}
\ifdim \wd\@tempboxa > 5.5in
\begin{center}
\begin{minipage}[t]{5.5in}
\addtolength{\baselineskip}{-0.95pt}
{\bf #1:} #2 \par
\addtolength{\baselineskip}{0.95pt}
\end{minipage}
\end{center}
\else
\hbox to\hsize{\hfil\box\@tempboxa\hfil}
\fi
\end{small}\par
}
\makeatother

\newcounter{oursection}

\newcounter{lecture}

\graphicspath{{img/}}

\usepackage{mathtools}
\usepackage{adjustbox}
\usepackage{booktabs}
\usepackage{csquotes}

\usepackage{hyperref}
\hypersetup{ colorlinks   = true,    % Colours links instead of ugly boxes
urlcolor     = blue,    % Colour for external hyperlinks
linkcolor    = black,    % Colour of internal links
citecolor    = blue      % Colour of citations
}

\usepackage{color}

\newcommand{\MS}[1]{\textbf{\color{cyan}{[MS: #1]}}}

\title{Automatic Generation of\\ Explicit Quadratic Programming Solvers}
\author{Maximilian Schaller 
\and Daniel Arnstr\"{o}m 
\and Alberto Bemporad
\and Stephen Boyd}
\date{June 12, 2025} 

\begin{document}
\maketitle

\begin{abstract}
We consider a family of convex quadratic programs in which 
the coefficients of the linear objective term and the righthand side of the
constraints are affine functions of a parameter.
It is well known that the solution of such a parametrized 
quadratic program is a piecewise affine function of the parameter.
The number of (polyhedral) regions in the solution map
can grow exponentially in problem size,
but when the number of regions is moderate, a so-called 
explicit solver is practical.
Such a solver computes the coefficients of the affine
functions and the linear inequalities defining the polyhedral regions 
offline; 
to solve a problem instance online it simply evaluates 
this explicit solution map.
Potential advantages of an explicit solver over a more general purpose
iterative solver can include transparency, interpretability, reliability,
and speed.
In this paper we describe how code generation
can be used to automatically generate an explicit solver from a 
high level description of a parametrized quadratic program.
Our method has been implemented in the open-source software 
CVXPYgen, which is part of CVXPY, a domain specific language for
general convex optimization.
\end{abstract}

\clearpage
\tableofcontents
\clearpage

\section{Introduction}
\subsection{Parametric convex optimization}
A parametric convex optimization problem can be written as
\BEQ\label{e-convex}
\begin{array}{ll}
\mbox{minimize} & f_0(x, \theta) \\
\mbox{subject to} & f_i(x, \theta) \le 0, \quad i=1,\ldots,m, \\
&h_i(x, \theta) = 0, \quad i=1,\ldots,q,
\end{array}
\EEQ
where $x \in \reals^n$ is the variable and $\theta \in \Theta \subseteq \reals^p$
is the parameter, \ie, data that is given and known whenever~\eqref{e-convex}
is solved.
The objective function $f_0$ and the inequality constraint functions
$f_i$, $i=1,\ldots,m$, are convex in $x$ and the equality constraint functions
$h_i$, $i=1,\ldots,q$, are affine in $x$, for any given value of
$\theta \in \Theta$~\cite{boyd2004convex}.
We refer to a solution of~\eqref{e-convex} as $x^\star(\theta)$ to emphasize its
dependence on the parameter $\theta$. Here we neglect that $x^\star$ might not
be unique or might not exist, and refer to the mapping from $\theta$ to $x$ 
as the solution map of the parametrized problem.

Convex optimization is used in various domains, including control systems
~\cite{borrelli2017predictive,rawlings2017model,kouvaritakis2016model,
keshavarz2014quadratic,wang2009fast,bemporad2007robust,boyd1991linear,
garcia1989model},
signal and image processing~\cite{chambolle2016introduction,
combettes2011proximal,mattingley2010real,zibulevsky2010l1},
and quantitative finance~\cite{palomar2025portfolio,boyd2024markowitz,
boyd2017multi,narang2013inside,grinold2000active,markowitz1952portfolio},
just to name a few that are particularly relevant for this work.

\paragraph{Explicit solvers for multiparametric programming.}
Traditionally, $x^\star(\theta)$ is evaluated using an iterative numerical method
that takes a given parameter value and computes an (almost) optimal point 
$x^\star(\theta)$~\cite{goulart2024clarabel,stellato2020osqp,odonoghue2016scs,domahidi2013ecos}.
We focus here on a very special case when $x^\star(\theta)$ can be 
expressed in closed form,
as an explicit function that maps a given value of $\theta$ directly to a 
solution $x^\star(\theta)$.
Such explicit solvers are practical for only some problems, and generally
only smaller instances, but when they are practical they can offer 
a number of advantages over generic iterative solvers.
Developing such solvers for parametric programs is now known as
multiparametric programming.

\subsection{Related work}

\paragraph{Multiparametric programming.}
Investigations of the theoretical properties of
parametric convex optimization problems date back to the 60s~\cite{MR64}
and were largely extended in the 80s~\cite{Fia83}.
After early work on explicitly solving parametric linear programs (LPs) in the context
of economics~\cite{gal1972multiparametric} in the 70s,
researchers started investigating the explicit solution of different convex 
optimization problems in the early 2000s.
Some of the first papers developed explicit solutions to
model predictive control problems based on 
quadratic programs (QPs)~\cite{bemporad2002explicit}
and LPs~\cite{bemporad2002model}, and general
algorithms were developed for explicitly solving QPs~\cite{bemporad2002explicit,tondel2003algorithm,PS10,GBN11}
and LPs~\cite{borrelli2003geometric}.
These are often cited as multiparametric linear or quadratic programming,
respectively, where the latter is often abbreviated as MPQP.
Further, people have worked on verifying the complexity of
such methods~\cite{cimini2017exact},
on approximate or suboptimal multiparametric 
programming~\cite{bemporad2003suboptimal, JG03, bemporad2006algorithm},
on multiparametric programming for linear complementary
problems~\cite{jones2006multiparametric},
and on synthesizing specialized hardware for multiparametric
programming~\cite{johansen2006hardware,BOPS11,ROLBHKS23}.

Software implementations include the
Hybrid Toolbox~\cite{HybTBX}, the Multi-Parametric Toolbox~\cite{MPT3},
the Model Predictive Control Toolbox~\cite{Bem15},
and the POP Toolbox~\cite{oberdieck2016pop} in Matlab,
the MPQP solver in the proprietary FORCES PRO software~\cite{domahidi2014forces},
and the PDAQP solver~\cite{arnstrom2024pdaqp} in Julia and Python,
which we use in this work.

Potential advantages of an explicit solver over a more general purpose
iterative solver can include transparency, interpretability, reliability, 
and speed.  Since the solver is essentially a lookup table, with
an explicit affine function associated with each region, there is
no question of convergence.  Indeed, we can explicitly determine the
maximum number of floating-point operations (FLOPS) required to compute
$x^\star(\theta)$ given $\theta$
\cite{cimini2017exact, arnstrom2024exact}. 
An explicit solver involves no division, so floating-point 
overflow or divide-by-zero 
exceptions cannot occur.  For the same reason, it is possible to store
the coefficients in a lower precision format such as 16-bit floating-point (FP16)
to reduce storage, and possibly to carry out the computations in 
lower precision as well, to increase speed or use low-cost electronic boards 
(at the cost of a modest decrease in accuracy).

The disadvantages of using an explicit solver
all relate to its worst-case exponential scaling with problem size.
This limits its practical use to relatively small problems,
which however do arise in many application areas.
Even when it is practical to use an explicit solver, the solver
data size can be large, since we must store the coefficients of the
explicit solution map.

\paragraph{Code generation for convex optimization.}
While typical convex optimization solvers are designed for general-purpose
computers~\cite{diamond2016cvxpy,goulart2024clarabel},
we are mostly interested in embedded applications with hard real-time constraints,
and also non-embedded applications where extreme speeds are required.

A \emph{code generator} heavily exploits the structure of the functions
in~\eqref{e-convex} and generates custom C code for solving the problem fast and
reliably for changing values of $\theta$, while fulfilling rules for safety-critical
code~\cite{holzmann2006power}.
Examples of code generators for iterative solvers
are CVXGEN~\cite{mattingley2012cvxgen} (for general QPs), CVXPYgen
~\cite{schaller2022embedded}, which interfaces with the OSQP
code generator~\cite{banjac2017embedded} and QOCOGEN~\cite{chari2025qoco}
(for QPs and second-order cone programs, respectively), and
acados~\cite{verschueren2021} and the proprietary
FORCES PRO~\cite{domahidi2014forces},
both specifically designed for QP-based and nonlinear control problems.
These code generators are used for many applications. CVXGEN, for example,
is used to guide and control all of the SpaceX first stage
landings \cite{blackmore2016autonomous}.

\paragraph{Domain-specific languages for optimization.}
Code generators typically accept problem specifications given
in a domain specific language (DSL) for convex optimization
\cite{lofberg2004yalmip,grant2014cvx}.
A DSL allows the user to describe the problem in a natural high level
human readable way, eliminating the effort (and risk of error) 
in transforming the problem to the standard form required by a solver.
For example there can be multiple variables or parameters,
with names that make sense in the application; the code
generator takes care of mapping these to our generic $x \in \reals^n$
and $\theta \in \Theta$.
Well-known DSLs include YALMIP~\cite{lofberg2004yalmip} and
CVX~\cite{grant2014cvx} in Matlab, CVXPY~\cite{diamond2016cvxpy} in Python,
Convex.jl~\cite{convexjl} and JuMP~\cite{Dunning2017jump} in Julia,
and CVXR~\cite{fu2017cvxr} in R.
We focus on CVXPY.

\subsection{Contribution}
In this paper, we adapt the code generator CVXPYgen
with the multiparametric explicit QP solver PDAQP
to generate explicit solvers for
convex optimization problems, described in CVXPY,
that can be reduced to QPs.
Along with C and C++ code for the generated solver, we generate a 
Python interface for rapid
prototyping and non-embedded applications.

We give four representative application examples, 
involving linear regression, 
power management in residential buildings, model predictive control, 
and financial portfolio optimization.
Our code generator accelerates the solve time for these 
examples by up to three orders of
magnitude compared to directly using CVXPY and its default QP solver,
with solve times down to hundreds of nanoseconds (on a standard laptop).

\subsection{Outline}

In \S\ref{s-explicit} we give a quick derivation of
the explicit solution map for a parametric QP.
In \S\ref{s-codegen} we explain how we generate source code
for an explicit solver that is described in the
modeling language CVXPY, and how to interface with the 
generated code.  We illustrate the explicit solver code generation
process in \S\ref{s-hello}.
In \S\ref{s-experiment} we assess the numerical performance
of the explicit solvers for several practical examples.
We report the time it takes to generate and compile the generated code,
the size of the resulting binary files, and the solve times.

\section{Explicit solution of parametric QPs}\label{s-explicit}

\subsection{Parametric QP}
We consider the QP
\BEQ\label{e-qp}
\begin{array}{ll}
\mbox{minimize} & (1/2) x^T P x + q ^T x \\
\mbox{subject to} & A x \leq b,
\end{array}
\EEQ
where $x \in \reals^n$ is the variable, and the data are
$P \in \symm^n_{++}$ 
(the set of symmetric positive definite $n \times n$ matrices),
$q\in \reals^n$, $A \in \reals^{m \times n}$, and $b \in \reals^m$.
The inequality in the constraints is elementwise.

Our focus is on the case when the data $P$ and $A$ are given, 
and $q$ and $b$ 
are affine functions of a parameter $\theta \in \reals^p$,
\BEQ\label{e-qb}
q = u + U \theta, \qquad
b = v + V \theta, 
\EEQ
where $u \in \reals^n$, $U\in \reals^{n \times p}$,
$v \in \reals^m$, and $V\in \reals^{m \times p}$ are given.
We refer to the QP \eqref{e-qp} parametrized by $\theta$ in 
\eqref{e-qb} as a parametrized QP.
Since the objective is strictly convex, there is at most one solution 
of the QP \eqref{e-qp} for each value of $\theta$.
We refer to the mapping from $\theta$ to the optimal $x$ (when it exists)
as the solution map of the parametrized QP \eqref{e-qp}.

\paragraph{Other QP forms.}
There are several other standard forms for a parametrized QP, 
both for analysis and as an interface to solvers
~\cite{bemporad2002explicit, tondel2003algorithm, arnstrom2024pdaqp},
but it is easy to translate between them by introducing additional
variables and constraints \cite{boyd2004convex}.

\paragraph{Constraints on the parameters.} In many applications
we are also given a set $\Theta \subseteq \reals^p$ of possible
parameter values.  For simplicity we ignore this, but occasionally
mention how this set of known possible values of $\theta$ can
be handled. When we do address the parameter set,
we assume it is a polyhedron.

\subsection{Optimality conditions}

\paragraph{Active constraints.}
We denote the $i$th row of $A$ as $a_i^T$,
so the inequality constraints
in \eqref{e-qp} can be expressed as $a_i^T x \leq b_i$,
$i=1, \ldots, m$.  We say that the $i$th inequality constraint is
tight or active if $a_i^T x = b_i$.
We let $\mathcal A = \{ i \mid a_i^T x = b_i \} 
\subseteq \{1, \ldots, m\}$ denote the set of active constraints
\cite{cimini2017exact,gill2019practical}
(which depends on $x$, $A$, and $b$).

Let $\lambda \in \reals_+^m$ denote a dual variable associated
with the linear inequality constraints in \eqref{e-qp}.
The optimality conditions for problem~\eqref{e-qp} are
\[
\begin{array}{l}
A x \leq b, \\
\lambda \geq 0,\\
Px + q + A^T \lambda = 0, \\
\lambda_i (a_i^T x - b_i) = 0, \quad i=1,\ldots,m.
\end{array}
\]
The first is primal feasibility; the second is nonnegativity of 
dual variables; the third is dual feasibility
(stationarity of the Lagrangian); and the last one is 
complementary slackness \cite[\S 5.5.3]{boyd2004convex}.

Let $\tilde A$, $\tilde b$, and $\tilde \lambda$ denote the 
row slices of $A$, $b$, and $\lambda$, respectively, corresponding
to the active constraints, \ie, $i \in \mathcal A$.
Let $\hat A$, $\hat b$, and $\hat \lambda$ denote the 
row slices of $A$, $b$, and $\lambda$, respectively, corresponding
to the inactive constraints, \ie, $i \not\in \mathcal A$.
By complementary slackness, we must have
$\lambda_i=0$ for $i \not\in \mathcal A$, so $\hat \lambda =0$.
With $\hat \lambda =0$, which we now assume, complementary slackness holds.
Since $\hat \lambda =0$, 
$A^T\lambda$ can be expressed as $\tilde A^T \tilde \lambda$, 
and dual feasibility can be expressed as
\[
Px+q + \tilde A^T \tilde \lambda = 0.
\]
Since $\mathcal A$ is the active set corresponding to $x$, we have
\[
\tilde A x = \tilde b.
\]
These two sets of linear equations can be summarized as the 
Karush-Kuhn-Tucker (KKT) system
\BEQ\label{e-kkt}
\left[\begin{array}{cc}
P & \tilde A^T \\
\tilde A & 0
\end{array}\right]
\left[\begin{array}{c}
x \\ \tilde \lambda
\end{array}\right]
=
\left[\begin{array}{c}
-q \\ \tilde b
\end{array}\right].
\EEQ
We assume that linear independence constraint qualification
(LICQ)~\cite{bertsekas1999nonlinear,nocedal2006numerical,boyd2004convex} holds,
\ie, that the rows of $\tilde A$ are linearly independent.
Then, it is well known that \eqref{e-kkt} can be uniquely solved
for $(x, \tilde \lambda)$~\cite[\S 10.1.1]{boyd2004convex},
\cite[\S 12.3]{boyd2018introduction}
and we re-write \eqref{e-kkt} as
\BEQ\label{e-kkt-solved}
\left[\begin{array}{c}
x \\ \tilde \lambda
\end{array}\right]
=
\left[\begin{array}{cc}
P & \tilde A^T \\
\tilde A & 0
\end{array}\right]^{-1}
\left[\begin{array}{c}
- q \\  \tilde b \end{array}\right].
\EEQ
This shows that knowledge of the active set $\mathcal A$ 
determines the primal and dual solutions of the QP \eqref{e-qp}.
We note that $(x,\tilde \lambda)$ are the solution of the 
problem of minimizing the objective of the QP subject to
the linear equality constraints $\tilde A x= \tilde b$.

From \eqref{e-kkt-solved} we see
that the solution is a linear function of the 
data $q$ and $b$, provided the active set does not change.
This implies that the solution is an affine function of the parameter
$\theta$, provided the active set does not change.

When $(x, \tilde \lambda)$ have the values \eqref{e-kkt-solved}
(with $\hat \lambda =0$), they satisfy complementary slackness 
and dual feasibility.  The remaining two optimality conditions,
primal feasibility and dual nonnegativity, can be expressed as 
\BEQ\label{e-polyhedron}
\left[\begin{array}{cc} \hat A & 0 \\ 0 &  -I \end{array}\right]
\left[\begin{array}{cc}
P & \tilde A^T \\
\tilde A & 0
\end{array}\right]^{-1}
\left[\begin{array}{c}
- q \\  \tilde b \end{array}\right]
\leq
\left[\begin{array}{c} \hat b \\ 0 \end{array}\right],
\EEQ
since $Ax \leq b$ holds (as equality) for rows
with $i \in \mathcal A$ and $\hat \lambda = 0$.
The inequality \eqref{e-polyhedron} is a set of linear 
inequalities in the data $b$ and $q$, and therefore
defines a polyhedron.
When $(b,q)$ is in this polyhedron, called the 
critical region associated with the active set $\mathcal A$, $(x,\lambda)$ given by 
the linear function \eqref{e-kkt-solved} are primal and dual optimal
for the QP \eqref{e-qp}. 

Since compositions of affine functions are affine, it follows that
the primal and dual solutions of the QP \eqref{e-qp} are (locally)
affine functions of $\theta$.  Since the inverse image of a polyhedron
under an affine mapping is a polyhedron, the values of $\theta$ 
over which this affine function gives the solution is also a 
polyhedron.
Thus the solution map is a piecewise affine function of 
$\theta$, with the polyhedral regions determined by the active set.
By the uniqueness of the solution, it is not difficult to prove that 
such a map is also continuous across region boundaries \cite{bemporad2002explicit}. 
This continuity property guarantees a certain degree of robustness to numerical errors 
when evaluating the map.

\subsection{Explicit parametric QP solver}
The optimality conditions discussed above suggest a naive explicit solver,
which can be practical when $m$ is small.
We search over all $2^m$ potential active sets. For each one
we compute $x$ and $\lambda$ via \eqref{e-kkt-solved}, and then check 
whether \eqref{e-polyhedron} holds.  If this happens, we have found 
the solution; if not, the problem is infeasible.

Now we consider the parameter dependence. For each of the $2^m$ 
potential active sets,
we can express \eqref{e-kkt-solved} as
\[
(x, \tilde \lambda) = F\theta +g, \qquad \hat \lambda = 0,
\]
where
\[
F =
\left[\begin{array}{cc}
P & \tilde A^T \\ \tilde A & 0
\end{array}\right]^{-1}
\left[\begin{array}{c}
-U \\ \tilde V
\end{array}\right], \qquad
g =
\left[\begin{array}{cc}
P & \tilde A^T \\ \tilde A & 0
\end{array}\right]^{-1}
\left[\begin{array}{c}
- u \\ \tilde v \end{array}\right],
\]
where $\tilde V$ and $\tilde v$ are the (row) slices of $V$
and $v$, respectively,
corresponding to $\mathcal A$.
We then express \eqref{e-polyhedron} in terms of $\theta$ as 
\[
H \theta \leq j,
\]
where
\[
H = \left[\begin{array}{cc} \hat A & 0 \\ 0 &  -I \end{array}\right]
F - \left[\begin{array}{c} \hat V \\ 0 \end{array}\right], \qquad
j= - \left[\begin{array}{cc} \hat A & 0 \\ 0 &  -I \end{array}\right]g
+ \left[\begin{array}{c} \hat v \\ 0 \end{array}\right] .
\]

For some choices of potential active sets, the inequalities $H\theta \leq j$
(together with $\theta \in \Theta$) are infeasible.
We drop these, and consider only the remaining $K$ sets of 
potential active sets, and label them as $k=1, \ldots, K$.
We can compute the coefficients of the piecewise affine solution map,
denoted $F_k$ and $g_k$, and their associated region, defined by
$H_k$ and $j_k$, before knowing the specific value of $\theta$.
Thus we have the explicit solution map
\BEQ\label{e-pwa-soln}
(x,\lambda) = F_k \theta + g_k \quad \mbox{when} \quad H_k\theta \leq j_k,
\quad k=1, \ldots, K.
\EEQ
Different existing MPQP solvers differ in how they avoid enumerating all $2^m$
active sets \cite[\S II-c]{arnstrom2024pdaqp}.
When $\theta \subseteq \Theta$ satisfies none of the inequalities above, the 
QP is infeasible.
The collection of coefficient matrices and vectors 
$F_k$, $g_k$, $H_k$, $j_k$, $k=1, \ldots, K$ gives an explicit 
representation of the solution map of the parametrized QP \eqref{e-qp}.
Since we can compute these matrices and vectors
(and determine the value of $K$) before we 
have specified $\theta$, we refer to computing these coefficients as the 
offline solve.
Evaluating \eqref{e-pwa-soln} for a given value of $\theta$ is called
the online solve.  Note that it involves no division.

The number of coefficients in the explicit solver is around
$K(n+m)p$, up to a factor of $K$ larger than the number of coefficients 
in the original problem, $n^2 + nm + (n+m)p$ (neglecting sparsity).
Even though $K$ can grow exponentially with $m$,
it is often practical for small
problems~\cite{bemporad2002explicit, tondel2003algorithm, cimini2017exact}.

\paragraph{Implementation.}
When we explicitly solve a parametrized QP in practice,
some of our initial assumptions can be relaxed. We can
handle positive semidefinite $P$ (instead of just positive definite);
we directly handle equality constraints; and we do not require LICQ.
(When $P$ is not positive definite, the solver provides \emph{a} solution,
rather than \emph{the} solution, since the solution need not be 
unique in this case.)

In the offline phase, implementations do not search all $2^m$ 
possible active sets, but instead find nonempty regions one by one.
This allows us to handle problems where $2^m$ is very large, but $K$,
the number of (nonempty) regions, is still moderate.

In the online solve, explicit solvers store 
the regions in a way that facilitates faster search
than a simple linear search over $k=1, \ldots, K$, typically involving
a pre-computed tree.
Other methods are used to either reduce the storage or increase the speed of
the online evaluations.

We use the specific solver PDAQP~\cite{arnstrom2024pdaqp},
which has several such accelerations implemented.
The offline solve is made more efficient by systematic searching over
neighboring regions.
The online solve benefits from a binary search tree for the 
search over regions~\cite{tondel2003evaluation}.
Complete details can be found in
\begin{center}
\url{https://github.com/darnstrom/pdaqp}.
\end{center}

\section{Code generation}\label{s-codegen}

\subsection{Domain-specific languages for optimization}

When solving a problem instance, DSLs perform a sequence of three steps.
First, the DSL transforms the user-defined problem into a form accepted by a standard
or canonical
solver. For example, a constraint like
\[
0 \leq x \leq 1
\]
for $x \in \reals$ is translated to
\[
A x \leq b, \quad 
A = \left[\begin{array}{r}
1 \\ -1
\end{array}\right], \quad
b = \left[\begin{array}{c}
1 \\ 0
\end{array}\right],
\]
as it appears in a canonical form like \eqref{e-qp}.
In a second step, a canonical solver (like PDAQP) is called to
solve the canonical problem. Ultimately, a solution for the user-defined
problem is retrieved from the solution returned by the canonical solver.

Typically, these three steps are performed every time a problem instance
is solved. We call such systems \emph{parser-solvers}.
When dealing with a parametric QP, whose structure does not change
between solves, repeated parsing, \ie, discovering how to reduce the
problem to canonical form, is unnecessary and usually inefficient.

\subsection{Code generation for explicitly solving QPs}
Consider an application where we solve many 
instances of a specific parametric QP,
possibly in an embedded system with hard real-time constraints.
For such applications, a \emph{code generator} makes more sense.

\begin{figure}
\centering
\includegraphics[width=0.8\columnwidth]{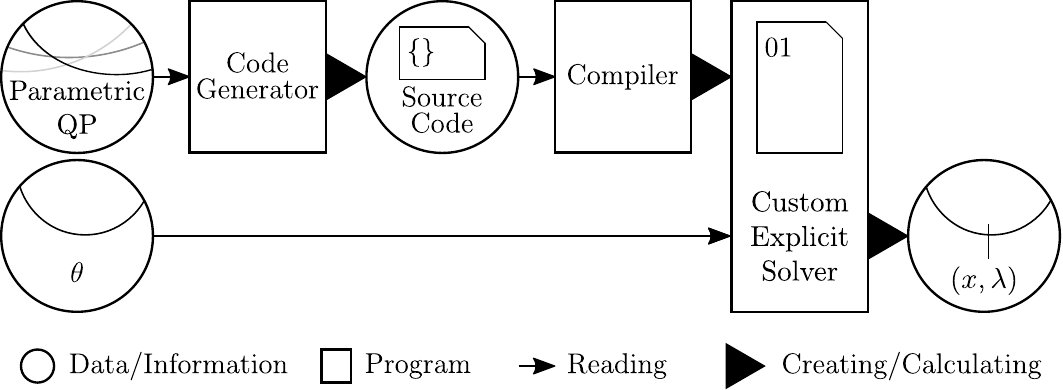}
\caption{Code generation for explicitly solving a parametric QP.}
\label{fig:codegen}
\end{figure}

As illustrated in figure~\ref{fig:codegen}, a 
code generator for explicitly solving QPs
takes as input the parametric QP, and generates
source code that is tailored for (explicitly) solving instances of that
parametric QP.
The source code is then compiled into an efficient custom solver,
which has a number of benefits compared to parser-solvers.
First, by exploiting the parametric structure and caching canonicalization,
the compiled solver becomes faster.
Second, the compiled solver can be deployed in embedded systems,
satisfying rules for safety-critical code~\cite{holzmann2006power}.

We extend the open-source code generator CVXPYgen~\cite{schaller2022embedded}
to generate code for explicitly solving QPs.
The QP is modeled with CVXPY before
CVXPYgen generates library-, allocation-, and division-free code for translating
between the user-defined problem and a canonical form (that of PDAQP in this case)
for explicitly solving QPs.
Open source code and full documentation for CVXPYgen and its 
explicit solve feature is available at
\begin{center}
\url{https://github.com/cvxgrp/cvxpygen}.
\end{center}

\paragraph{Disciplined parametrized programming.}
We require the QP to be modeled in CVXPY using disciplined parametrized programming (DPP)
~\cite{agrawal2019differentiable}. The DPP rules mildly restrict how parameters may enter
the objective and constraints. Generally speaking, if parameters enter all expressions in an affine
way, then the problem is DPP-compliant. Details on the DPP rules can be found at
\url{https://www.cvxpy.org}.

\paragraph{Canonicalization.}
When a problem is modeled according to DPP, parameter canonicalization and
solution retrieval are affine mappings,
\[
\theta = C  \theta^\text{user} + c, \qquad
x^\text{user} = R x + r,
\]
where $\theta^\text{user}$ and $x^\text{user}$ are the user-defined parameter
and variable, respectively.
The matrices $C$ and $R$ are typically very sparse.
The retrieval matrix $R$ is usually wide, \ie,
there are more canonical variables than user-defined variables, since
the canonicalization step introduces auxiliary variables~\cite{agrawal2018rewriting}.
CVXPYgen generates code for
the respective sparse matrix-vector multiplications. In fact, the solution retrieval
$Rx + r$ often reduces to simple pointing to memory in C.

\paragraph{Parameter constraints.}
The user may specify constraints on parameters.
Suppose that the values
of the user-defined parameter $\theta^\text{user}$ lie
in the set $\Theta^\text{user}$.
Then, similarly to how the parameter $\theta^\text{user}$ 
itself is canonicalized,
the parameter constraints are translated to the canonical set of
possible parameters $\Theta$.

\paragraph{Limitations.}
In PDAQP the number of inequality constraints $m$ is limited to 1024 so 
the regions can be efficiently represented as a bit string.
If $K$, or the size of the data in the explicit solver, exceed a given
limit, the offline phase is terminated with a warning.

\section{Hello world} \label{s-hello}
Here we present a simple example to illustrate how explicit solver 
code generation works.
Consider the parametric QP
\BEQ\label{e-helloworld}
\begin{array}{ll}
\mbox{minimize} & \|X \beta + v \ones - y\|_2^2 \\
\mbox{subject to} & \beta \geq 0,
\end{array}
\EEQ
where $\beta \in \reals^d$ and $v \in \reals$ are the variables,
$\theta^\text{user} = y \in \reals^p$
is the parameter with $\Theta^\text{user} = \{ y \mid l \leq y \leq u \}$,
and $\ones$ denotes the vector with all entries one.
The bounds $l \in \reals^p$ and $u \in \reals^p$ and the
matrix $X \in \reals^{p \times d}$ are given data.

\subsection{Modeling and code generation}
The problem can be formulated in CVXPY as shown in
figure~\ref{code-modeling}, up to line~11. Note that in line~10,
we specify $\Theta^\text{user}$ with standard CVXPY constraints.
We generate the explicit solver in line~13.
\begin{figure}
\lstset{language=Python,
numbers=left,
xleftmargin=0.12\columnwidth,
linewidth=0.88\columnwidth}
\begin{lstlisting}[frame=lines]
import cvxpy as cp
from cvxpygen import cpg

d, p, X, l, u = ...
beta = cp.Variable(d, name='beta')
v = cp.Variable(name='v')
y = cp.Parameter(p, name='y')

obj = cp.Minimize(cp.sum_squares(X @ beta + v - y))
constr = [beta >= 0, l <= y, y <= u]
prob = cp.Problem(obj, constr)

cpg.generate_code(prob, solver='explicit')
\end{lstlisting}
\caption{Generating an explicit solver with CVXPYgen.}
\label{code-modeling}
\end{figure}

\subsection{C interface}
Figure~\ref{code-C} shows how the generated explicit solver can be used in~C.
In line~7, the first entry of the parameter
$y$ is updated to the value $1.2$. The function \verb|cpg_update_y|
ensures that the parameter values are within their pre-specified
limits (otherwise, it maps the given value back onto $\Theta$).
In line~8, the problem is solved explicitly with the \verb|cpg_solve|
function. In lines~9 and~10, respectively,
the first entry of the optimal coefficients $\beta$ and the optimal $v$ 
is read and printed. In line~11, the resulting objective value is
calculated and printed. Note that the calculation of the objective value
is kept in a separate function from the solve function,
for maximal efficiency.
\begin{figure}
\lstset{language=C,
numbers=left,
xleftmargin=0.18\columnwidth,
linewidth=0.82\columnwidth}
\begin{lstlisting}[frame=lines]
#include <stdio.h>
#include "cpg_workspace.h"
#include "cpg_solve.h"

int main(int argc, char *argv[]){
	
	cpg_update_y(0, 1.2);
	cpg_solve();
	printf("%f\n", CPG_Result.prim->beta[0]);
	printf("%f\n", CPG_Result.prim->v);
	printf("%f\n", cpg_obj());
	
	return 0;
	
}
\end{lstlisting}
\caption{Using the explicit solver.}
\label{code-C}
\end{figure}

Figure~\ref{code-C-structs} shows the C structs
that store the result. In the example above, we access the
primal solution \verb|beta| and \verb|v| via the \verb|prim| field
of the result struct \verb|CPG_Result|.
The dual variables in \verb|CPG_Dual_t| are named according
to the index in the list of CVXPY constraints.
\begin{figure}
\lstset{language=C,
numbers=left,
xleftmargin=0.15\columnwidth,
linewidth=0.85\columnwidth}
\begin{lstlisting}[frame=lines]
typedef struct {
	cpg_float  *beta;    // primal variable beta
	cpg_float   v;       // primal variable v
} CPG_Prim_t;

typedef struct {
	cpg_float  *d0;      // dual variable d0
} CPG_Dual_t;

typedef struct {
	CPG_Prim_t *prim;    // primal solution
	CPG_Dual_t *dual;    // dual solution
} CPG_Result_t;
\end{lstlisting}
\caption{Data structure of explicit solver result.}
\label{code-C-structs}
\end{figure}

\subsection{CVXPY interface}
We consider a small instance of the parametric QP~\eqref{e-helloworld}
with $d=2$, $p=3$, the entries of $X$ generated IID from $\mathcal{N}(0, 1)$,
$l = 0$, and $u = \ones$.
We assign the entries of $y$ randomly between $0$ and $1$ and
solve the problem three times: with CVXPY using the iterative OSQP
solver~\cite{stellato2020osqp}, with CVXPYgen using OSQP,
and with CVXPYgen using the explicit PDAQP solver.

Figure~\ref{code-solving} shows the comparison, demonstrating
how the explicit solver can be used via its auto-generated CVXPY interface.
Starting in line 15, we show that the primal and dual solutions
and the objective values are all close, respectively.
\begin{figure}
\lstset{language=Python,
numbers=left,
xleftmargin=0.15\columnwidth,
linewidth=0.85\columnwidth}
\begin{lstlisting}[frame=lines]
from code_osqp.cpg_solver import cpg_solve
prob.register_solve('gen_OSQP', cpg_solve)

from code_explicit.cpg_solver import cpg_solve
prob.register_solve('gen_explicit', cpg_solve)

def print_result():
  print(f'v:    {v.value}')
  print(f'beta: {beta.value}')
  print(f'dual: {constr[0].dual_value}')
  print(f'obj:  {obj.value}')

y.value = [0.6, 0.8, 0.2]

prob.solve(solver='OSQP')
print_result()

# v:    0.916741
# beta: [0.000000 0.288547]
# dual: [0.750370 0.000000]
# obj:  0.059628

prob.solve(method='gen_OSQP')
print_result()

# v:    0.916741
# beta: [0.000000 0.288547]
# dual: [0.750370 0.000000]
# obj:  0.059628

prob.solve(method='gen_explicit')
print_result()

# v:    0.916740
# beta: [0.000000 0.288546]
# dual: [0.750370 0.000000]
# obj:  0.059628
\end{lstlisting}
\caption{Using the explicit solver in CVXPY.}
\label{code-solving}
\end{figure}

\clearpage
\section{Applications}\label{s-experiment}
In this section we report timing and code size details for some typical
application examples.
In each case we compare CVXPYgen using the explicit
PDAQP solver to CVXPYgen using the iterative 
OSQP solver and standard CVXPY using OSQP.
When using CVXPYgen with the OSQP solver, we use OSQP's code
generation feature, which caches the factorization of the KKT
system, for accelerated solving~\cite{banjac2017embedded}.
When using OSQP in CVXPY or CVXPYgen, we set both the relative
and absolute tolerances to $10^{-4}$ (the default in CVXPY).
We run the experiments on an Apple M1 Pro, compiling with Clang
at optimization level 3.
For iterative and explicit code generation, we report solve times in C,
and the overall time when solving from Python via the auto-generated
CVXPY interface.
We also give the time it takes to generate and compile the code.

\subsection{Monotone regression}

We consider the monotone regression problem
\[
\begin{array}{ll}
\mbox{minimize} & \|A x - b\|_2^2 \\
\mbox{subject to} & x_1 \leq x_2 \leq \cdots \leq x_d,
\end{array}
\]
where $x \in \reals^d$ is the variable and $\theta^\text{user} = b \in \reals^q$
is the parameter, with $\Theta^\text{user} = [-1,1]^q$.
The matrix $A \in \reals^{q \times d}$ is given.
(This is called monotone regression since the components
$x_i$ are constrained to be monotonically nondecreasing.)

\paragraph{Problem instances.}
We consider $d = 5$ and $q = 10$.
The entries of $A$ are generated IID from $\mathcal{N}(0, 1)$,
and we generate 100 problem instances where $b$ is sampled
uniformly from $\Theta^\text{user}$.

\paragraph{Results.}
We obtain $n=15$ variables, $m=d-1=4$
inequality constraints, and $p=q=10$ parameters.
(In this case, the canonicalization step introduced $n-d=10$ auxiliary variables.)
We find $K=16$ regions (which equals $2^m$, the maximum possible number of
active sets).
Table~\ref{tab:simple_qp} shows the average solve times and binary sizes.
\begin{table}
\begin{minipage}{\textwidth}
\centering
\begin{tabular}{l|rr|rr|r}
& Solve (Python) & Solve (C) & Gen. + compile & Gen. & Binary size \\ \hline
CVXPY & 0.6089 ms & -- & -- & -- & -- \\
CVXPYgen OSQP & 0.1764 ms & 0.1257 ms & 5.7 s & 0.1 s &  80 KB  \\
CVXPYgen explicit & 0.0127 ms & 0.0004 ms & 13.7 s & 9.5 s & 15 KB  \\
\hline
\end{tabular}
\caption{Timing and binary sizes for monotone regression problem.}
\label{tab:simple_qp}
\end{minipage}
\end{table}

\subsection{Power management}

\paragraph{Problem.}
A nonnegative electric power load $L$ is served by a 
PV (photovoltaic solar panel) system, 
a storage battery, and a grid connection~\cite{byrne2017energy, nnorom2025aging}.
We denote the solar power as $s$,
the battery power as $b$, and the grid power as $g$.
These three power sources supply the load, so we have
\[
L = s + b + g.
\]
The PV power satisfies $0 \leq s \leq S$, where $S \geq 0$
is the available PV power.
The battery power satisfies $-C \leq b \leq D$, where $D>0$
is the maximum possible discharge power and $C>0$ is the maximum
possible charge power.
The grid power satisfies $g \geq 0$, \ie, we cannot sell 
power back to the grid.
The (positive) price of the grid power is $P$, so the grid 
cost is $Pgh$, where $h$ is the duration of one time period,
over which we hold the power values constant.

The battery state of charge at the beginning of the time period 
is denoted $q$, and satisfies
$0 \leq q \leq Q$, where $Q$ is the battery capacity.
At the beginning of the next time period the battery charge 
is $q^+ = q - h b$.  We must have $0 \leq q^+ \leq Q$.

We take the cost function
\[
Pgh + \alpha (q^+-q^\text{tar})^2 + \beta b^2,
\]
where $\alpha$ and $\beta$ are given and positive,
and $q^\text{tar}$ is a given target battery charge value.

To choose the powers we solve the QP
\[
\begin{array}{ll}
\mbox{minimize} & Pgh + \alpha (q^+ - q^\text{tar})^2 + \beta b^2\\
\mbox{subject to} & L = s + b + g, \\
& 0 \leq s \leq S, \quad  -C \leq b \leq D, \quad g \geq 0, \\
& q^+ = q - hb, \quad 0 \leq q^+ \leq Q,
\end{array}
\]
where $s$, $b$, $g$, and $q^+$ are the variables, and
$\theta^\text{user} = (L, S, P, q)$ are parameters. The remaining constants,
$C$, $D$, $h$, $Q$, $q^\text{tar}$, $\alpha$, and $\beta$,
are known.
We take
\[
\Theta^\text{user} = [0, 1] \times [0, 0.5] \times [1, 2] \times [0, Q].
\]

\paragraph{Problem instances.}
We set  $C=D=1$, $h=0.05$, $Q=1$, $q^\text{tar}=0.5$, and
$\alpha = \beta = 0.1$.
We generate 100 problem instances where $(L, S, P, q)$ is sampled
uniformly from $\Theta^\text{user}$.

\paragraph{Results.}
After canonicalization, there are $n=5$ variables (including one
auxiliary variable), $m=7$
inequality constraints, and $p=4$ parameters. We find $K=5$ regions.
Table~\ref{tab:energy} shows the average solve times and binary sizes.
\begin{table}
\begin{minipage}{\textwidth}
\centering
\begin{tabular}{l|rr|rr|r}
& Solve (Python) & Solve (C) & Gen. + compile & Gen. & Binary size \\ \hline
CVXPY & 0.6564 ms & -- & --  & -- & -- \\
CVXPYgen OSQP &0.1261 ms & 0.0541 ms & 5.1 s & 0.1 s & 75 KB \\
CVXPYgen explicit & 0.0207 ms & 0.0001 ms & 12.6 s & 8.7 s & 10 KB \\
\hline
\end{tabular}
\caption{Timing and binary sizes for the power management problem.}
\label{tab:energy}
\end{minipage}
\end{table}

\subsection{Model predictive control}
We consider the linear dynamical system~\cite{boyd1991linear}
\[
z_{t+1} = A z_t + B u_t, \quad t=0,1,\ldots,H-1,
\]
where $z_t \in \reals^{n_z}$ is the state and 
$u_t \in \reals^{n_u}$ is the input, which must satisfy 
$\|u_t\|_\infty \leq 1$.
The matrices $A \in \reals^{n_z \times n_z}$ and 
$B \in \reals^{n_z \times n_u}$ are given.
We solve the model predictive control problem~\cite{garcia1989model,kouvaritakis2016model}
\[
\begin{array}{ll}
\mbox{minimize} & z_H^T P z_H + \sum_{t=0}^{H-1} \left(z_t^T Q z_t + u_t^T R u_t\right)  \\
\mbox{subject to} & z_{t+1} = A z_t + B u_t, \quad t=0,\ldots, H-1, \\
& \|u_t\|_\infty \leq 1, \quad  t=0,\ldots, H-1, \\
& z_0 = z^\text{init},
\end{array}
\]
where $z_0, \ldots, z_H$ and $u_0, \ldots, u_{H-1}$ are the variables
and $\theta^\text{user} = z^\text{init}$ is the parameter.
We take $\Theta^\text{user} = [-1,1]^{n_z}$.
The objective matrices $P \in \symm^{n_z}_{++}$, $Q \in \symm^{n_z}_{++}$,
and $R \in \symm^{n_u}_{++}$ (along with $A$ and $B$) are given.

\paragraph{Problem instances.}
We consider $n_z=6$ states, $n_u=1$ input, and a horizon length of $H = 5$.
We construct $A$ by sampling its diagonal entries from $\mathcal{N}(0, 1)$
and its off-diagonal entries IID from $\mathcal{N}(0, 0.01)$,
before scaling the whole matrix such that $A$ has spectral radius $1$.
The entries of the input matrix $B$ are sampled IID from $\mathcal{N}(0, 0.001)$.
We set the controller weights to $Q = I$ and $R = 0.1 I$, and compute $P$
as the solution to the algebraic Riccati equation
associated with the infinite-horizon problem~\cite{kwakernaak1972linear}.
We generate 100 problem instances where
the entries of $z^\text{init}$ are generated uniformly from $\Theta^\text{user}$.

\paragraph{Results.}
After canonicalization, we have $n=77$ variables (of which $36$ % = (H+1)n_z
are auxiliary variables), $m=10$ inequality
constraints, and $p = n_z = 6$ parameters.
We find $K=63$ regions.
Table~\ref{tab:control} shows the average solve times and binary sizes.
\begin{table}
\begin{minipage}{\textwidth}
\centering
\begin{tabular}{l|rr|rr|r}
& Solve (Python) & Solve (C) & Gen. + compile & Gen. & Binary size \\ \hline
CVXPY & 1.102 ms & -- & -- & -- & --\\
CVXPYgen OSQP &0.875 ms & 0.790 ms & 5.2 s & 0.1 s & 110 KB \\
CVXPYgen explicit & 0.025 ms & 0.001 ms & 13.7 s & 9.1 s & 93 KB \\
\hline
\end{tabular}
\caption{Timing and binary sizes for the model predictive control problem.}
\label{tab:control}
\end{minipage}
\end{table}

\subsection{Portfolio optimization}\label{s-portfolio}
We construct a financial portfolio consisting of holdings in $N$ assets
~\cite{grinold2000active,narang2013inside,palomar2025portfolio}.
We represent the holdings relative to the total (positive) portfolio value,
in terms of nonnegative weights $w \in \reals_+^N$, where $\ones^T w = 1$,
with $w_i$ being the fraction of the (positive) total portfolio value 
invested in asset $i$.
With estimated mean annualized asset returns $\mu \in \reals^N$
\cite{grinold2000active}, 
the estimated mean annualized portfolio return is $\mu^T w$.
The variance or risk of the portfolio return is $w^T \Sigma w$,
where $\Sigma \in \symm^N_{++}$ is an estimate for the covariance
matrix of the annualized asset returns.
Our objective is to maximize the risk-adjusted expected annualized return
\[
\mu^T w - \gamma w^T \Sigma w,
\]
where $\gamma > 0$ is the risk-aversion factor.
To find the portfolio we solve the Markowitz
problem~\cite{markowitz1952portfolio,boyd2024markowitz}
\[
\begin{array}{ll}
\mbox{maximize} & \mu^T w - \gamma w^T \Sigma w \\
\mbox{subject to} & \ones^T w = 1, \quad w \geq 0, 
\end{array}
\]
where the portfolio weights $w \in \reals^N$
are the variable and the parameter is $\theta^\text{user} = \mu$ with
$\Theta^\text{user} = [-1, 1]^N$.
This means that we expect no annualized returns beyond $\pm 100 \%$.
The covariance matrix $\Sigma \in \symm^N_{++}$ and the
risk-aversion factor $\gamma > 0$ are given.

\paragraph{Problem instances.}
We take $N = 7$ assets. 
To obtain data we choose the $7$ stocks with
the largest market capitalization as of January 1, 2017, listed
in table~\ref{tab:stocks}.
\begin{table}
\begin{minipage}{\textwidth}
\centering
\begin{tabular}{ll}
Ticker symbol & Company \\ \hline
AAPL & Apple \\
AMZN & Amazon \\
BRK.A & Berkshire Hathaway \\
FB (now META) & Facebook (now Meta) \\
GOOGL & Alphabet \\
%JNJ & Johnson \& Johnson \\
%JPM & JPMorgen Chase \\
MSFT & Microsoft \\
%TCEHY & Tencent \\
XOM & ExxonMobil \\
\hline
\end{tabular}
\caption{Stocks used in the portfolio optimization problem.}
\label{tab:stocks}
\end{minipage}
\end{table}
We compute $\Sigma$ by first taking the sample covariance of the $7$ assets'
daily returns in the years 2017 and 2018, and then annualizing the result.
We choose $\gamma = 2$.
We generate 250 problem instances where $\mu$ is taken as the
one-year trailing average of returns for 250 trading days in the year 2019.

\paragraph{Results.}
We have $n=N=7$ variables, $m=N=7$ inequality
constraints, and $p =N= 7$ parameters.
We find $K=127$ regions, only one less than the maximum $2^m$ potential
combinations of investing in an asset or not, since the constraint
$\ones^T w = 1$ prevents $w = 0$ (not investing in any asset).
Table~\ref{tab:portfolio} shows the average solve times and binary sizes.
\begin{table}
\begin{minipage}{\textwidth}
\centering
\begin{tabular}{l|rr|rr|r}
& Solve (Python) & Solve (C) & Gen. + compile & Gen. & Binary size \\ \hline
CVXPY & 0.5441 ms & -- & -- & -- & -- \\
CVXPYgen OSQP & 0.0502 ms & 0.0070 ms & 5.1 s & 0.1 s & 76 KB \\
CVXPYgen explicit & 0.0113 ms & 0.0005 ms & 20.8 s & 16.5 s & 234 KB \\
\hline
\end{tabular}
\caption{Timing and binary sizes for the portfolio optimization problem.}
\label{tab:portfolio}
\end{minipage}
\end{table}

\clearpage
\section{Conclusions}
We have added new functionality to the code generator CVXPYgen
that generates an explicit solver (in C) for a parametrized convex 
optimization problem, when that is tractable.
The user can prototype a problem in CVXPY, with code close to the math
and convenient names for multiple variables and parameters, using
a generic iterative solver; a change of one option in code generation
will generate an explicit solver for the parametrized problem.
For typical (small) problems from various application domains,
our numerical experiments show the generated explicit solvers 
exhibit solve times at (or below) one microsecond, giving up to three
orders of magnitude speedup over an iterative solver.

{%\small
\bibliography{refs}

\newcommand{\etalchar}[1]{$^{#1}$}
\begin{thebibliography}{RMD{\etalchar{+}}17}

\bibitem[AA24]{arnstrom2024pdaqp}
D.~Arnström and D.~Axehill.
\newblock A high-performant multi-parametric quadratic programming solver.
\newblock In {\em 2024 IEEE 63rd Conference on Decision and Control (CDC)},
  pages 303--308, 2024.

\bibitem[AAB{\etalchar{+}}19]{agrawal2019differentiable}
A.~Agrawal, B.~Amos, S.~Barratt, S.~Boyd, S.~Diamond, and Z.~Kolter.
\newblock Differentiable convex optimization layers.
\newblock In {\em Advances in Neural Information Processing Systems (NeurIPS)},
  volume~32, 2019.

\bibitem[ABA24]{arnstrom2024exact}
D.~Arnstr{\"o}m, D.~Broman, and D.~Axehill.
\newblock Exact worst-case execution-time analysis for implicit model
  predictive control.
\newblock {\em IEEE Transactions on Automatic Control}, 69(10):7190--7196,
  2024.

\bibitem[AVDB18]{agrawal2018rewriting}
A.~Agrawal, R.~Verschueren, S.~Diamond, and S.~Boyd.
\newblock A rewriting system for convex optimization problems.
\newblock {\em Journal of Control and Decision}, 5(1):42--60, 2018.

\bibitem[BB91]{boyd1991linear}
S.~Boyd and C.~Barratt.
\newblock {\em Linear controller design: {L}imits of performance}, volume~78.
\newblock Citeseer, 1991.

\bibitem[BBD{\etalchar{+}}17]{boyd2017multi}
S.~Boyd, E.~Busseti, S.~Diamond, R.~Kahn, K.~Koh, P.~Nystrup, and J.~Speth.
\newblock Multi-period trading via convex optimization.
\newblock {\em Foundations and Trends in Optimization}, 3(1):1--76, 2017.

\bibitem[BBM{\etalchar{+}}02]{bemporad2002model}
A.~Bemporad, F.~Borrelli, M.~Morari, et~al.
\newblock Model predictive control based on linear programming -- the explicit
  solution.
\newblock {\em IEEE transactions on automatic control}, 47(12):1974--1985,
  2002.

\bibitem[BBM03]{borrelli2003geometric}
F.~Borrelli, A.~Bemporad, and M.~Morari.
\newblock Geometric algorithm for multiparametric linear programming.
\newblock {\em Journal of optimization theory and applications}, 118:515--540,
  2003.

\bibitem[BBM17]{borrelli2017predictive}
F.~Borrelli, A.~Bemporad, and M.~Morari.
\newblock {\em Predictive control for linear and hybrid systems}.
\newblock Cambridge University Press, 2017.

\bibitem[Bem04]{HybTBX}
A.~Bemporad.
\newblock {Hybrid Toolbox - User's Guide}, 2004.
\newblock \url{http://cse.lab.imtlucca.it/~bemporad/hybrid/toolbox}.

\bibitem[Bem15]{Bem15}
A.~Bemporad.
\newblock A multiparametric quadratic programming algorithm with polyhedral
  computations based on nonnegative least squares.
\newblock {\em IEEE Transactions on Automatic Control}, 60(11):2892--2903,
  2015.

\bibitem[Ber99]{bertsekas1999nonlinear}
D.~Bertsekas.
\newblock {\em Nonlinear Programming}.
\newblock Athena Scientific, 2nd edition, 1999.

\bibitem[BF03]{bemporad2003suboptimal}
A.~Bemporad and C.~Filippi.
\newblock Suboptimal explicit receding horizon control via approximate
  multiparametric quadratic programming.
\newblock {\em Journal of optimization theory and applications}, 117:9--38,
  2003.

\bibitem[BF06]{bemporad2006algorithm}
A.~Bemporad and C.~Filippi.
\newblock An algorithm for approximate multiparametric convex programming.
\newblock {\em Computational optimization and applications}, 35:87--108, 2006.

\bibitem[BJK{\etalchar{+}}24]{boyd2024markowitz}
S.~Boyd, K.~Johansson, R.~Kahn, P.~Schiele, and T.~Schmelzer.
\newblock Markowitz portfolio construction at seventy.
\newblock {\em Journal of Portfolio Management}, 50(8):117--160, 2024.
\newblock Also available at \url{https://arxiv.org/pdf/2401.05080}.

\bibitem[Bla16]{blackmore2016autonomous}
L.~Blackmore.
\newblock Autonomous precision landing of space rockets.
\newblock In {\em Frontiers of Engineering: {R}eports on Leading-Edge
  Engineering from the 2016 Symposium}, volume~46, pages 15--20. The Bridge
  Washington, DC, 2016.

\bibitem[BM07]{bemporad2007robust}
A.~Bemporad and M.~Morari.
\newblock Robust model predictive control: A survey.
\newblock In {\em Robustness in identification and control}, pages 207--226.
  Springer, 2007.

\bibitem[BMDP02]{bemporad2002explicit}
A.~Bemporad, M.~Morari, V.~Dua, and E.~Pistikopoulos.
\newblock The explicit linear quadratic regulator for constrained systems.
\newblock {\em Automatica}, 38(1):3--20, 2002.

\bibitem[BNC{\etalchar{+}}17]{byrne2017energy}
R.~Byrne, T.~Nguyen, D.~Copp, B.~Chalamala, and I.~Gyuk.
\newblock Energy management and optimization methods for grid energy storage
  systems.
\newblock {\em IEEE Access}, 6:13231--13260, 2017.

\bibitem[BOPS11]{BOPS11}
A.~Bemporad, A.~Oliveri, T.~Poggi, and M.~Storace.
\newblock Ultra-fast stabilizing model predictive control via canonical
  piecewise affine approximations.
\newblock {\em IEEE Transactions on Automatic Control}, 56(12):2883--2897,
  2011.

\bibitem[BSM{\etalchar{+}}17]{banjac2017embedded}
G.~Banjac, B.~Stellato, N.~Moehle, P.~Goulart, A.~Bemporad, and S.~Boyd.
\newblock Embedded code generation using the {OSQP} solver.
\newblock In {\em IEEE Conference on Decision and Control}, pages 1906--1911.
  IEEE, 2017.

\bibitem[BV04]{boyd2004convex}
S.~Boyd and L.~Vandenberghe.
\newblock {\em Convex Optimization}.
\newblock Cambridge University Press, 2004.

\bibitem[BV18]{boyd2018introduction}
S.~Boyd and L.~Vandenberghe.
\newblock {\em Introduction to Applied Linear Algebra: Vectors, Matrices, and
  Least Squares}.
\newblock Cambridge University Press, 2018.

\bibitem[CA25]{chari2025qoco}
G.~Chari and B.~A{\c{c}}{\i}kme{\c{s}}e.
\newblock {QOCO}: {A} quadratic objective conic optimizer with custom solver
  generation.
\newblock {\em arXiv preprint arXiv:2503.12658}, 2025.

\bibitem[CB17]{cimini2017exact}
G.~Cimini and A.~Bemporad.
\newblock Exact complexity certification of active-set methods for quadratic
  programming.
\newblock {\em IEEE Transactions on Automatic Control}, 62(12):6094--6109,
  2017.

\bibitem[CP11]{combettes2011proximal}
P.~Combettes and J.~Pesquet.
\newblock Proximal splitting methods in signal processing.
\newblock {\em Fixed-point algorithms for inverse problems in science and
  engineering}, pages 185--212, 2011.

\bibitem[CP16]{chambolle2016introduction}
A.~Chambolle and T.~Pock.
\newblock An introduction to continuous optimization for imaging.
\newblock {\em Acta Numerica}, 25:161--319, 2016.

\bibitem[DB16]{diamond2016cvxpy}
S.~Diamond and S.~Boyd.
\newblock {CVXPY}: {A} {P}ython-embedded modeling language for convex
  optimization.
\newblock {\em Journal of Machine Learning Research}, 17(83):1--5, 2016.

\bibitem[DCB13]{domahidi2013ecos}
A.~Domahidi, E.~Chu, and S.~Boyd.
\newblock {ECOS}: {A}n {SOCP} solver for embedded systems.
\newblock In {\em European Control Conference (ECC)}, pages 3071--3076. IEEE,
  2013.

\bibitem[DHL17]{Dunning2017jump}
I.~Dunning, J.~Huchette, and M.~Lubin.
\newblock {JuMP}: {A} modeling language for mathematical optimization.
\newblock {\em SIAM Review}, 59(2):295--320, 2017.

\bibitem[DJ14]{domahidi2014forces}
A.~Domahidi and J.~Jerez.
\newblock {FORCES} {P}rofessional. {E}mbotech {G}mb{H}, 2014.

\bibitem[Fia83]{Fia83}
A.~Fiacco.
\newblock {\em Introduction to Sensitivity and Stability Analysis in Nonlinear
  Programming}.
\newblock Academic Press, London, U.K., 1983.

\bibitem[FNB20]{fu2017cvxr}
A.~Fu, B.~Narasimhan, and S.~Boyd.
\newblock {CVXR}: {A}n {R} package for disciplined convex optimization.
\newblock {\em Journal of Statistical Software}, 94(14):1--34, 2020.

\bibitem[GB14]{grant2014cvx}
M.~Grant and S.~Boyd.
\newblock {CVX}: {M}atlab software for disciplined convex programming, version
  2.1, 2014.

\bibitem[GBN11]{GBN11}
A.~Gupta, S.~Bhartiya, and P.~Nataraj.
\newblock A novel approach to multiparametric quadratic programming.
\newblock {\em Automatica}, 47(9):2112--2117, 2011.

\bibitem[GC24]{goulart2024clarabel}
P.~Goulart and Y.~Chen.
\newblock Clarabel: An interior-point solver for conic programs with quadratic
  objectives.
\newblock {\em arXiv preprint arXiv:2405.12762}, 2024.

\bibitem[GK00]{grinold2000active}
R.~Grinold and R.~Kahn.
\newblock {\em Active portfolio management}.
\newblock McGraw Hill New York, 2000.

\bibitem[GMW19]{gill2019practical}
P.~Gill, W.~Murray, and M.~Wright.
\newblock {\em Practical optimization}.
\newblock SIAM, 2019.

\bibitem[GN72]{gal1972multiparametric}
T.~Gal and J.~Nedoma.
\newblock Multiparametric linear programming.
\newblock {\em Management Science}, 18(7):406--422, 1972.

\bibitem[GPM89]{garcia1989model}
C.~Garcia, D.~Prett, and M.~Morari.
\newblock Model predictive control: {T}heory and practice -- a survey.
\newblock {\em Automatica}, 25(3):335--348, 1989.

\bibitem[HKJM13]{MPT3}
M.~Herceg, M.~Kvasnica, C.~N. Jones, and M.~Morari.
\newblock Multi-parametric toolbox 3.0.
\newblock {\em European Control Conference (ECC)}, pages 502--510, 2013.

\bibitem[Hol06]{holzmann2006power}
G.~Holzmann.
\newblock The power of 10: {R}ules for developing safety-critical code.
\newblock {\em Computer}, 39(6):95--99, 2006.

\bibitem[JG03]{JG03}
T.~Johansen and A.~Grancharova.
\newblock Approximate explicit constrained linear model predictive control via
  orthogonal search tree.
\newblock {\em IEEE Transactions on Automatic Control}, 58(5):810--815, 2003.

\bibitem[JJST06]{johansen2006hardware}
T.~Johansen, W.~Jackson, R.~Schreiber, and P.~Tondel.
\newblock Hardware synthesis of explicit model predictive controllers.
\newblock {\em IEEE Transactions on control systems technology},
  15(1):191--197, 2006.

\bibitem[JM06]{jones2006multiparametric}
C.~Jones and M.~Morrari.
\newblock Multiparametric linear complementarity problems.
\newblock In {\em Proceedings of the 45th IEEE Conference on Decision and
  Control}, pages 5687--5692. IEEE, 2006.

\bibitem[KB14]{keshavarz2014quadratic}
A.~Keshavarz and S.~Boyd.
\newblock Quadratic approximate dynamic programming for input-affine systems.
\newblock {\em International Journal of Robust and Nonlinear Control},
  24(3):432--449, 2014.

\bibitem[KC16]{kouvaritakis2016model}
B.~Kouvaritakis and M.~Cannon.
\newblock {\em Model Predictive Control}.
\newblock Springer, 2016.

\bibitem[KS72]{kwakernaak1972linear}
H.~Kwakernaak and R.~Sivan.
\newblock {\em Linear optimal control systems}.
\newblock Wiley-InterScience New York, 1972.

\bibitem[L\"04]{lofberg2004yalmip}
J.~L\"ofberg.
\newblock {YALMIP}: {A} toolbox for modeling and optimization in {Matlab}.
\newblock In {\em IEEE International Conference on Robotics and Automation
  (ICRA)}, pages 284--289. IEEE, 2004.

\bibitem[Mar52]{markowitz1952portfolio}
H.~Markowitz.
\newblock Portfolio selection.
\newblock {\em Journal of Finance}, 7(1):77--91, 1952.

\bibitem[MB10]{mattingley2010real}
J.~Mattingley and S.~Boyd.
\newblock Real-time convex optimization in signal processing.
\newblock {\em IEEE Signal Processing Magazine}, 27(3):50--61, 2010.

\bibitem[MB12]{mattingley2012cvxgen}
J.~Mattingley and S.~Boyd.
\newblock {CVXGEN}: {A} code generator for embedded convex optimization.
\newblock {\em Optimization and Engineering}, 13:1--27, 2012.

\bibitem[MR64]{MR64}
O.~Mangasarian and J.~Rosen.
\newblock Inequalities for stochastic nonlinear programming problems.
\newblock {\em Operations Research}, 12:143--154, 1964.

\bibitem[Nar13]{narang2013inside}
R.~Narang.
\newblock {\em Inside the Black Box: A Simple Guide to Quantitative and
  High-frequency Trading}.
\newblock John Wiley \& Sons, 2013.

\bibitem[NOBL25]{nnorom2025aging}
O.~Nnorom, G.~Ogut, S.~Boyd, and P.~Levis.
\newblock Aging-aware battery control via convex optimization.
\newblock {\em arXiv preprint arXiv:2505.09030}, 2025.

\bibitem[NW06]{nocedal2006numerical}
J.~Nocedal and S.~Wright.
\newblock {\em Numerical Optimization}.
\newblock Springer, 2nd edition, 2006.

\bibitem[OCPB16]{odonoghue2016scs}
B.~O'Donoghue, E.~Chu, N.~Parikh, and S.~Boyd.
\newblock Conic optimization via operator splitting and homogeneous self-dual
  embedding.
\newblock {\em Journal of Optimization Theory and Applications},
  169(3):1042--1068, 2016.

\bibitem[ODP{\etalchar{+}}16]{oberdieck2016pop}
R.~Oberdieck, N.~Diangelakis, M.~Papathanasiou, I.~Nascu, and E.~Pistikopoulos.
\newblock {POP} -- {P}arametric optimization toolbox.
\newblock {\em Industrial \& Engineering Chemistry Research},
  55(33):8979--8991, 2016.

\bibitem[Pal25]{palomar2025portfolio}
D.~Palomar.
\newblock {\em Portfolio Optimization: {T}heory and Application}.
\newblock Cambridge University Press, 2025.

\bibitem[PS10]{PS10}
P.~Patrinos and H.~Sarimveis.
\newblock A new algorithm for solving convex parametric quadratic programs
  based on graphical derivatives of solution mappings.
\newblock {\em Automatica}, 46(9):1405--1418, 2010.

\bibitem[RMD{\etalchar{+}}17]{rawlings2017model}
J.~Rawlings, D.~Mayne, M.~Diehl, et~al.
\newblock {\em Model Predictive Control: Theory, Computation, and Design}.
\newblock Nob Hill Publishing Madison, WI, 2017.

\bibitem[ROL{\etalchar{+}}23]{ROLBHKS23}
A.~Ravera, A.~Oliveri, M.~Lodi, A.~Bemporad, W.~Heemels, E.~Kerrigan, and
  M.~Storace.
\newblock Co-design of a controller and its digital implementation: {T}he
  {MOBY-DIC2} toolbox for embedded model predictive control.
\newblock {\em IEEE Transactions on Control Systems Technology},
  31(6):2871--2878, 2023.

\bibitem[SBD{\etalchar{+}}22]{schaller2022embedded}
M.~Schaller, G.~Banjac, S.~Diamond, A.~Agrawal, B.~Stellato, and S.~Boyd.
\newblock Embedded code generation with {CVXPY}.
\newblock {\em IEEE Control Systems Letters}, 6:2653--2658, 2022.

\bibitem[SBG{\etalchar{+}}20]{stellato2020osqp}
B.~Stellato, G.~Banjac, P.~Goulart, A.~Bemporad, and S.~Boyd.
\newblock {OSQP}: {A}n operator splitting solver for quadratic programs.
\newblock {\em Mathematical Programming Computation}, 12(4):637--672, 2020.

\bibitem[TJB03a]{tondel2003algorithm}
P.~T{\o}ndel, T.~Johansen, and A.~Bemporad.
\newblock An algorithm for multi-parametric quadratic programming and explicit
  {MPC} solutions.
\newblock {\em Automatica}, 39(3):489--497, 2003.

\bibitem[TJB03b]{tondel2003evaluation}
P.~T{\o}ndel, T.~Johansen, and A.~Bemporad.
\newblock Evaluation of piecewise affine control via binary search tree.
\newblock {\em Automatica}, 39(5):945--950, 2003.

\bibitem[UMZ{\etalchar{+}}14]{convexjl}
M.~Udell, K.~Mohan, D.~Zeng, J.~Hong, S.~Diamond, and S.~Boyd.
\newblock Convex optimization in {J}ulia.
\newblock In {\em 2014 first workshop for high performance technical computing
  in dynamic languages}, pages 18--28. IEEE, 2014.

\bibitem[VFK{\etalchar{+}}21]{verschueren2021}
R.~Verschueren, G.~Frison, D.~Kouzoupis, J.~Frey, N.~van Duijkeren, A.~Zanelli,
  B.~Novoselnik, T.~Albin, R.~Quirynen, and M.~Diehl.
\newblock acados -- a modular open-source framework for fast embedded optimal
  control.
\newblock {\em Mathematical Programming Computation}, pages 1--37, 2021.

\bibitem[WB09]{wang2009fast}
Y.~Wang and S.~Boyd.
\newblock Fast model predictive control using online optimization.
\newblock {\em IEEE Transactions on Control Systems Technology},
  18(2):267--278, 2009.

\bibitem[ZE10]{zibulevsky2010l1}
M.~Zibulevsky and M.~Elad.
\newblock {$L_1$-$L_2$} optimization in signal and image processing.
\newblock {\em IEEE Signal Processing Magazine}, 27(3):76--88, 2010.

\end{thebibliography}
}

\end{document}